\documentclass[secthm]{elsart}

\usepackage{epsfig}

\usepackage{amsmath,multirow,amssymb}

\theoremstyle{definition}
\newtheorem{eg}[thm]{Example}
\renewcommand{\qed}{\hfill $\blacksquare$}
\allowdisplaybreaks

\newlength{\leftparboxindent}

\def \constA {c_1}
\def \constB {c_2}
\def \constC {c_3}
\def \constD {c_4}
\def \constE {c_5}
\def \constF {c_6}
\def \constG {c_7}
\def \constH {c_8}
\def \constI {c_9}
\def \constJ {c_{10}}
\def \constK {c_{11}}
\def \constL {c_{12}}
\def \constM {c_{13}}
\def \constN {c_{14}}

\begin{document}

\setlength{\leftparboxindent}{\textwidth}

\addtolength{\leftparboxindent}{-.5in}


\addtolength{\oddsidemargin}{-.25in}\addtolength{\evensidemargin}{-.25in}
\addtolength{\marginparwidth}{.55in}\addtolength{\marginparsep}{-.1in}

\begin{frontmatter}

\title{Two-batch liar games on a general bounded channel}

\author[Ellis]{Robert B. Ellis\corauthref{Ecorauth}}
\ead{rellis@math.iit.edu}
\author[Nyman]{Kathryn L. Nyman}

\corauth[Ecorauth]{Author to whom Correspondence should be sent.}
\ead{knyman@math.luc.edu}

\address[Ellis]{Department of Applied Mathematics, Illinois Institute of
Technology, Chicago, IL 60616}
\address[Nyman]{Department of Mathematics and Statistics, Loyola University
Chicago, Chicago, IL 60626}

\begin{abstract}
We consider an extension of the 2-person R\'enyi-Ulam liar game in
which lies are governed by a channel $C$, a set of allowable lie
strings of maximum length $k$. Carole selects $x\in[n]$, and Paul
makes $t$-ary queries to uniquely determine $x$. In each of $q$
rounds, Paul weakly partitions $[n]=A_0\cup \cdots \cup A_{t-1}$
and asks for $a$ such that $x\in A_a$. Carole responds with some
$b$, and if $a\neq b$, then $x$ accumulates a lie $(a,b)$.
Carole's string of lies for $x$ must be in the channel $C$.  Paul
wins if he determines $x$ within $q$ rounds. We further restrict
Paul to ask his questions in two off-line batches. We show that
for a range of sizes of the second batch, the maximum size of the
search space $[n]$ for which Paul can guarantee finding the
distinguished element is $\sim t^{q+k}/(E_k(C)\binom{q}{k})$ as
$q\rightarrow\infty$, where $E_k(C)$ is the number of lie strings
in $C$ of maximum length $k$.  This generalizes previous work of
Dumitriu and Spencer, and of Ahlswede, Cicalese, and Deppe. We
extend Paul's strategy to solve also the pathological liar
variant, in a unified manner which gives the existence of
asymptotically perfect two-batch adaptive codes for the channel
$C$.
\end{abstract}

\end{frontmatter}

\section{Introduction}

We consider a generalization of the R\'enyi-Ulam liar game,
originating in \cite{R61} and \cite{U76}. In this 2-player ``20
questions'' game, Paul may ask 20 Yes-No questions in order to
identify a distinguished element $x$ from a set
$[n]:=\{1,\ldots,n\}$, where Carole answers ``Yes'' or ``No'' and
is allowed to lie at most once. If Paul is allowed $q$ questions,
he can identify $x$ provided $n \lesssim 2^q/(q+1)$ (see
\cite{P87}). Restricting Carole to always tell the truth reduces
the game to binary search.  An equivalent coding theoretic
formulation is block coding over a noisy binary symmetric channel
with noiseless feedback \cite{B68}.

The basic R\'enyi-Ulam liar game has these parameters: search
space $[n]$, number of questions $q$, and Carole's maximum number
of lies $k$.
In \cite{DS05b}, Dumitriu and Spencer determined the first-term
asymptotics of the following extension: instead of binary Yes-No
questions, Paul asks $t$-ary questions, and Carole has a set of
lie types (e.g., ``Yes'' when the truth is ``No'') from which she
may draw up to $k$ times with repetition.
Furthermore, Paul asks his questions in two batches, receiving
Carole's answers at the end of a batch.  In \cite{ACD08},
Ahlswede, Cicalese, and Deppe extended this result to weighted
lies, with bounded total weight.
In \cite{EY04}, the first author and Yan
introduced the pathological variant of the liar game, in which
Paul wins provided at least one element in the search space
survives being eliminated, with Carole playing adversely.

In this paper, we simultaneously unify and extend
\cite{DS05b,ACD08,EY04} as follows. Our channel is a finite list
of strings of lies of varying type, from which Carole selects one
string to apply its lies in order and interspersed among her
responses. Every candidate $y\in[n]$ has a game lie string
generated by Carole's $q$ responses from the perspective of $y$
being the element Paul seeks; if $y$'s string is not in the list,
then $y$ is eliminated from the search space.
Furthermore, Paul is constrained to the aforementioned
two-batch question strategy. We solve asymptotically both the
original and the pathological variants for the optimal $n$ for
which Paul can win in terms of $q$, giving unified winning
strategies for Paul in the original and pathological variants
which correspond to asymptotically perfect adaptive codes.
This, our main result, is given as Theorems
\ref{thm:PaulWinCond}-\ref{thm:CaroleWinCond} in Section
\ref{sec:mainResult}, with proofs deferred to Section
\ref{sec:proof}, followed by concluding remarks in Section
\ref{sec:conclusion}.
A list of principal notation appears in Table
\ref{tab:principalNotation} after Section \ref{sec:defn}, and
the beginning of Section \ref{sec:PaulWinCond} is a technical
outline of Paul's unified winning strategies.

Our general channel condition is natural because it encompasses
the previously studied binary and $t$-ary liar games on a
symmetric, asymmetric, or weighted channel. It also specializes to
the binary unidirectional channel, in which lies may be of one
type or the other but not mixed (see Example \ref{exa:unidir});
for a bounded number of lies in arbitrary position, we believe
that the most general previous result is
\cite{ACD08}. Furthermore, the pathological variant appears to be
less studied, with results only in the binary asymmetric and
binary symmetric cases \cite{EY04,EPY05}.
Finally, requiring Paul to ask questions in two batches with a
range of possible sizes for the second batch provides intuition
about the fully nonadaptive one-batch case, which includes
$k$-error-correcting codes and radius $k$ covering codes, in the
original and pathological variants, respectively.

\section{Definitions and main result\label{sec:defn}}

The R\'enyi-Ulam liar game is a 2-player perfect information game
with integer parameters $n,q,k \geq 0$ and $t\geq 2$, played over
a $t$-ary communication channel $C$ of order $k$, which we now
define. The {\em lies} for the alphabet $T:=\{0,\ldots,t-1\}$ are
the set
$$\mathcal{L}(t):=\{(a,b)\in T\times T:a\neq b\}.$$
A {\em lie string} is a finite ordered list of lies, that is, an
element of the language $\mathcal{L}(t)^*:=\cup_{i\geq
0}\mathcal{L}(t)^i$.
For our
purposes, a {\em $t$-ary communication channel of order $k$} is an
arbitrary subset
$$
    C \subseteq \bigcup_{i= 0}^k\mathcal{L}(t)^i,
$$
such that $C\cap\mathcal{L}(t)^k\neq \emptyset$, representing
the usable lie strings of the game.  We denote the order of $C$
by $o(C)$. The unique element and {\em empty string} $\epsilon$
of $\mathcal{L}(t)^0$ may or may not be in $C$. The {\em
length} $|u|$ of $u\in\mathcal{L}(t)^*$ is simply the number of
lies in $u$. The {\em concatenation} of
$u,v\in\mathcal{L}(t)^*$ is defined as
$uv:=(a_1,b_1)\cdots(a_{j'},b_{j'})$ when
$u=(a_1,b_1)\cdots(a_{j},b_{j})$ and
$v=(a_{j+1},b_{j+1})\cdots(a_{j'},b_{j'})$.

Paul and Carole play a $q$-round game on the set
$[n]:=\{1,\ldots,n\}$.
Each $y\in[n]$ begins the game with lie string $\epsilon$. To
start each round, Paul weakly partitions $[n]$ into $t$ parts by
choosing a {\em question} $(A_0,\ldots,A_{t-1})$ such that
$[n]=A_0\dot\cup\cdots\dot\cup A_{t-1}$, where $\dot\cup$ denotes
disjoint union. Carole completes the round by responding with her
{\em answer}, an index $j\in T$. If $y\in A_j$, then $y$
accumulates no lie.  Otherwise $y\in A_i$ for some $i\neq j$, and
$(i,j)$ is post-pended to $y$'s current lie string.
The {\em game response string} is the ordered sequence $w'_1\cdots
w'_q\in T^q$ of Carole's responses. The {\em game lie string} for
$y$ is its final accumulated sequence of lies
$(a_1,b_1)\cdots(a_j,b_j)\in\cup_{i=0}^q\mathcal{L}(t)^i$.
By definition, the $b_i$'s must appear, in order, as a
substring of
the game response string.
\begin{equation}
\begin{array}{r@{\hspace{.2in}}ccccccccc}
\mathrm{Truthful~response~string~for~}y\mathrm{:}& w_1 & \cdots
    & w_{i_1} & \cdots & w_{i_\ell}
    & \cdots & w_{i_j} & \cdots & w_q\\
\multirow{2}{*}{Game~lie~string~for~$y$:} & & & a_1 & \cdots
    & a_{\ell} & \cdots & a_j & & \\
& & & b_1 & \cdots
    & b_{\ell} & \cdots & b_j & & \\
\mathrm{Game~response~string:} & w'_1 & \cdots & w'_{i_1} & \cdots
    & w'_{i_\ell} & \cdots & w'_{i_j} &
\cdots & w'_q
\end{array} \label{eqn:applyLieString}
\end{equation}
Here, $w_{i_\ell} = a_\ell$ and $w'_{i_\ell}=b_\ell$ for all $1\leq
\ell\leq j$, and $w_i=w'_i$ for all other indices.
If $y$'s game lie string is in $C$, then $y$ {\em survives} the
game; otherwise $y$ is {\em disqualified} ({\em eliminated}). More
broadly, at any
given round, 
$y$ survives iff its final lie string might still be in $C$.
Rather than requiring Carole to choose $x$ at the beginning, we
equivalently allow her to update her choice of $x$, lie string,
and game response string at any time.
Paul wins the {\em original variant} iff after $q$ rounds at most
one element (candidate for $x$) survives, and he wins the {\em
pathological variant} iff after $q$ rounds at least one element
survives. For the second variant, we think of a capricious Carole
lying ``pathologically'' in order to eliminate all elements as
quickly as possible.  Carole plays an adversarial strategy in both
variants and wins if Paul does not. In a {\em fully adaptive}
game, Paul receives Carole's answer each round before forming his
next question.  We will usually restrict Paul to a {\em two-batch
strategy} consisting of $q_1$ questions in the first batch and
$q_2$ questions in the second batch.  Carole responds to an entire
batch at once after receiving all questions in the batch.
\begin{defn}
The {\em $(n,q_1,q_2,C)$-game} is the two-batch original liar game
with search space $[n]$ on a $t$-ary channel $C$ of order
$o(C)<\infty$, 
with $q_1$ and $q_2$ questions in the first and
second batches, respectively.  The {\em $(n,q_1,q_2,C)^*$-game} is
the two-batch pathological liar game with the same parameters.
\end{defn}

For the binary channel $C=\{(0,1)\}$, Carole must lie since
$\epsilon\notin C$. But Paul may win the original variant
regardless of $[n]$ by setting $A_0=\emptyset$ and $A_1=[n]$ every
round. To avoid such trivial winning strategies, we constrain $C$
as follows.
\begin{defn}[Non-degenerate channel]\label{def:nondegen}
The channel $C$ is {\em non-degenerate} with respect to the original
variant provided either \vspace{-.1in}
\begin{enumerate}
\item $\epsilon\in C$, or \item for all $a\in T$, there exists a
lie string $u\in C$ with $u=(a,b_1)\cdots (a,b_j)$;
\end{enumerate} \vspace{-.1in}
and  is {\em non-degenerate} with respect to the
pathological variant provided either \vspace{-.1in}
\begin{enumerate}
\item $\epsilon\in C$, or \item for all $b\in T$, there exists a
lie string $u\in C$ with $u=(a_1,b)\cdots (a_j,b)$.
\end{enumerate}
\end{defn}
In the above example, $C=\{(0,1)\}$ had no $u$ of the form
$(1,b_1)\cdots(1,b_j)$. In the pathological variant, unless $C$ is
non-degenerate there exists $b\in T$ with no $u\in C$ of the form
$(a_1,b)\cdots(a_j,b)$. Carole wins regardless of $[n]$ by always
answering $b$, thereby eliminating every $y\in[n]$.
The fully adaptive case needs a more careful definition of
non-degeneracy (though $\epsilon\in C$ suffices), which we leave
to future work.

\subsection{The main result\label{sec:mainResult}}

For a $t$-ary channel $C$ of order $k$ and for $0\leq j\leq k$,
define the number
$$
    E_j(C) := \left|C\cap \mathcal{L}(t)^j\right|
$$
of length $j$ lie strings in $C$. Our main result is that, for
$q_2$ sufficiently bounded, the asymptotic optimal $n$ for
which Paul can win the $(n,q_1,q_2,C)$-game or
$(n,q_1,q_2,C)^*$-game depends on $E_k(C)$ and not on $C$
itself.
For convenience, we separate the main result into bounds for Paul
and  bounds for Carole.

\begin{thm} \label{thm:PaulWinCond}
Let $C$ be an order $k$ channel, let $f(q)$ be nonnegative with
$f(q)\rightarrow\infty$, and let $q_1+q_2=q$.
There exist constants $\constA, \constB$
such that if $(\ln q)^{3/2}f(q)\leq q_2\leq \constA
q^{k/(2k-1)}$ and
\begin{equation}\label{eqn:PaulOrigWinCond}
n\leq \frac{t^{q+k}}{E_k(C)\binom{q}{k}}\left(1-\constB\frac{\sqrt{\ln
    q}}
    {q_2^{1/3}}
    \right),
\end{equation}
then, for $q$ large enough, Paul can win the $(n,q_1,q_2,C)$-game.
If, in addition, $C$ is non-degenerate with respect to the
pathological variant, then there exists a constant $\constC$
such that if $(\ln q)^{3/2}f(q)\leq q_2\leq \constA
q^{k/(2k-1)}$ and
\begin{equation}\label{eqn:PaulWinCond}
n\geq
\frac{t^{q+k}}{E_k(C)\binom{q}{k}}\left(1+\constC\frac{\sqrt{\ln
    q}}{q^{1/3}}\right),
\end{equation}
then, for $q$ large enough, Paul can win the $(n,q_1,q_2,C)^*$-game.
\end{thm}
\begin{thm} \label{thm:CaroleWinCond}
Let $C$ be an order $k$ channel, let $f(q)$ be nonnegative with
$f(q)\rightarrow\infty$, and let $q_1+q_2=q$ with
$\max(q_1,q_2)\gg\min(q_1,q_2)\geq(\ln q)^{3/2} f(q)$.
There exist constants $\constD,\constE$ such that if $C$ is
non-degenerate with respect to the original variant and
\begin{align}
n \geq \frac{t^{q+k}}{E_k(C)\binom qk}
    \left(1 + \frac{\constD\min(q_1,q_2)}{q} + \constE
    \frac{\sqrt{\ln q}}{\max(q_1,q_2)^{1/3}}
    \right), \nonumber
\end{align}
then, for $q$ sufficiently large, Carole can win the
$(q_1,q_2,n,C)$-game. Regardless of the choice of $q_1$ and $q_2$,
there exists a constant $\constF$ such that if
\begin{align}
n \leq \frac{t^{q+k}}{E_k(C)\binom qk}
    \left ( 1 -\frac{\constF\sqrt{\ln q}}{q^{1/3}} \right ), \nonumber
\end{align}
then, for $q$ sufficiently large, Carole can win the $(q_1,q_2, n,
C)^*$-game.
\end{thm}

The above constants depend on $k$ but not on $q$, $q_1$, or
$q_2$.
We defer proofs until Section \ref{sec:proof}.  The proof of
Theorem \ref{thm:PaulWinCond} builds on that of Theorem 1.2 of
\cite{DS05b}, which gives the optimal $n$, up to the first
asymptotic term in $q$, for which Paul can win the original
game variant when $C=\cup_{j=0}^k \mathcal{L}'(t)^j$, for some
fixed $\mathcal{L}'(t)\subseteq \mathcal{L}(t)$.  We borrow
their two-batch structure, but extend it to handle a more
general channel, the pathological game variant, and a wider
range of second batch size $q_2$. Most original is our unified
treatment in the key Theorem \ref{thm:PaulPackVol} of Paul's
winning strategies in both variants, which proves the existence
of asymptotically perfect adaptive codes for any non-degenerate
channel $C$. These codes correspond to packings within
coverings of the $t$-ary hypercube $T^q$ of the sets of game
response strings for which individual elements of the search
space survive (like Hamming balls for nonadaptive codes). Our
proof for when Carole has a winning strategy borrows from
\cite{DS05b} but applies the two-batch structure in the
original variant as is necessary to be consistent with the
definition of a non-degenerate channel. A motivation for our
definition of $C$ was the following example.
\begin{eg}\label{exa:unidir}
In a liar game over a binary asymmetric ($Z$-) channel, Carole may
lie with ``Yes'' when the correct answer is ``No'' but not
vice-versa. In the companion asymmetric channel, only a ``No'' to
``Yes'' lie is allowed.  The 2-lie unidirectional channel
$C=\{\epsilon, (0,1),$ $(0,1),$ $(0,1)(0,1),$ $(1,0)$ $(1,0)\}$
may be interpreted as Paul knowing that the game is being played
over one of the asymmetric channels with $k=2$, but not which. In
prior work on the fully nonadaptive case (e.g.,
\cite{FvT92,OS06}), $C$ is called the ``unidirectional error''
channel.
\end{eg}
Substantially more general channels are possible.  For example, by
setting $C = \{(0,1),(1,0), (0,1)(1,0), (1,0)(0,1)\}$, we force
Carole to lie (as $\epsilon\notin C$), and require that if she
lies twice, her second lie must be of the opposite type.

\subsection{Suffix channels and the game state vector}

At any given round in the game, an element $y\in[n]$ has
accumulated a partial game lie string $u\in\mathcal{L}(t)^*$.  We
define the {\em suffix channel} $S_C(u)$ to be the set of all ways
to extend to a game lie string in the channel.  Formally,
$$
 S_C(u) :=
 \{v:
uv\in C\}, \qquad \mbox{and} \qquad \mathcal{S}(C) := \{S_C(u):u\in
\mathcal{L}(t)^*\}
$$
is the set of suffix channels of $C$.
Disqualified elements $y\in[n]$ have suffix channel $\emptyset$.
We track each $y\in [n]$ via its suffix channel at any given
round.

\begin{defn}\label{def:stateVector}
The {\em game state vector} after a given round in the original or
pathological game is the vector $(x_{C'}:C'\in\mathcal{S}(C))$
indexed by the suffix channels of $C$, where $x_{C'}$ counts the
number of elements of $[n]$ whose accumulated lie string $u$
satisfies $S_C(u)=C'$. The {\em coarsened state vector},
$(x_0,x_1,\ldots,x_k)$, is obtained from the game state vector by
grouping elements with suffix channels of the same order, so that
$$
    x_i := \mathop{\sum_{C'\in \mathcal{S}(C)}}_{o(C')=k-i}
    x_{C'} \qquad \mbox {for $0 \leq i \leq k$}.
$$
\end{defn}

At the start of our game, $x_C=n$ and $x_{C'}=0$ for $C'\neq C$.
If, after $q'\leq q$ rounds, $y\in[n]$ has partial game lie string
$u$, $y$ survives the entire game iff $y$ survives the game on the
last $q-q'$ rounds with respect to the suffix channel $S_C(u)$.
The element $y$ has been disqualified after $q'\leq q$ rounds iff
$S_C(u)=\emptyset$ or $q-q'<\min\{|v|:v\in S_C(u)\}$,
that is, if there are no strategies of questions by Paul and
answers by Carole in the last $q-q'$ rounds in which the
partial lie string of $y$ can be completed to obtain a game lie
string in $C$.
Thus we have the following.
\begin{lem}\label{lem:middleWinCond}
Let $0\leq q_1\leq q$ and $q_1+q_2=q$.  Given that the state
vector is $(x_{C'}:C'\in\mathcal{S}(C))$ after $q_1$ rounds,  Paul
wins the entire game, in either variant, iff he wins the
$q_2$-round game starting with state vector
$(x_{C'}:C'\in\mathcal{S}(C))$. \hfill \qed
\end{lem}
The state vector after $q_1$ rounds is a snapshot of the game
regardless of adaptive or two-batch questioning. Setting $q_1=q$,
Paul wins the original (pathological) variant iff after $q$ rounds
$\sum_{C',\epsilon\in C'}x_{C'}\leq 1$ ($\geq 1$), as the empty
lie string must be in an element's suffix channel for it to
survive with no questions left.

We conclude with channel statistics needed in Theorem
\ref{thm:PaulPackVol}.  The number of prefixes $u$ of an order $i$
suffix channel of $C$ is
$p_i(C):=|\{u\in\mathcal{L}(t)^*:o(S_C(u))=i\}|$.  Refining
$p_i(C)$ by the length of $u$, set $p_{i}^{(j)}(C):=|\{u \in
\mathcal{L}(t)^j:o(S_C(u))=i\}|$. Note that $p_{i}^{(j)}(C)=0$
when $j>k-i$, and that $p_i(C)=\sum_{j=0}^{k-i}p_{i}^{(j)}(C)$.

\begin{table}
\caption{Principal notation \label{tab:principalNotation}}
\begin{center}
\begin{tabular}[h]{ll}
\hline
$[n]$, $y$ & search space $[n]:=\{1,\ldots,n\}$ and typical element
    $y\in[n]$ \\
$q$, $q_1$, $q_2$ & number of total, first batch, and second batch
    questions, resp. \\
$(A_0,\ldots,A_{t-1})$ & $t$-ary question by Paul weakly partitioning
    $[n]$, where $t\geq 2$ \\
$T^Q$ & response strings $\{0,\ldots,t-1\}^Q$ for a batch of $Q$ questions\\
$(a,b)\in T\times T$ & w.r.t.~some $y$, truth $a$ and response $b$; a
    ``lie'' when $a\neq b$\\
$\mathcal{L}(t)$, $\mathcal{L}(t)^i$ & set of all lies of $T\times T$,
  and of all lie strings of length $i$ \\
$C\subseteq \cup_{i=0}^k\mathcal{L}(t)^i$ &  $t$-ary channel of order
     $o(C)=k$, with $C\cap \mathcal{L}(t)^k\neq \emptyset$\\
$u\in \mathcal{L}(t)^*$ & usually, the accumulated lie
    string of some $y\in[n]$ \\
$uv\in C$ & lie string of surviving $y$; $u$, $v$ from first, second
    batches, resp. \\
$w$, $w'$ & truthful response string for some $y$, actual response
    string \\
$z$, $z'$ & same as previous but usually for second batch \\
$S_C(u)$ & suffix channel of $u$, the set of $v$ with $uv\in C$ a legal lie
    string \\
$x_{C'}$, $x_i$ & counts $y\in [n]$ surviving with suffix channel
    $C'$, $o(C')=k-i$ \\
$(M,r)$-balanced & all $w\in T^Q$ with $\leq \frac{1}{t}
    \lceil Q/M\rceil +r$ per section of each letter of $T$ \\
$T^Q(M,r)$ & set of $(M,r)$-balanced strings of $T^Q$ \\
$\{w':w\stackrel{u}{\rightarrow}w'\}$ & set of all $w'$ arising from $w$
    under application of $u$ as in \eqref{eqn:applyLieString} \\
$B(w,C)$ & $C$-shadow with stem $w$; union over $u\in C$ of
    $\{w':w\stackrel{u}{\rightarrow}w'\}$ \\
$(x_{C'}:C'\in \mathcal{S}(C))$ & intra-batch state vector indexed by
    all suffix channels $C'$\\
$(x_i:0\leq i\leq k)$ & previous vector additively grouped by $i=k-o(C')$\\
$\alpha$ & number of $y$ Paul assigns to each
    $T^{q_1}(M_1,r_1)$, original variant\\
$\alpha'$ & additional number of $y$ for previous, pathological variant\\
$A_t(Q,2j+1)$ & maximum size of a $t$-ary length $Q$ code with packing radius $j$\\
\hline
\end{tabular}
\end{center}
\end{table}

\section{Proof of the main result\label{sec:proof}}

We begin with a notion of balanced strings of $T^q$ that have
nearly equal frequencies of each letter of $T$.  All game response
strings for which a typical $y\in[n]$ survives will be balanced.
For a 1-batch game, this set is a $C$-shadow, which is generalized
in Section \ref{sec:codingTheoretic} from a Hamming ball. In
Section \ref{sec:PaulWinCond}, Paul's winning strategies for both
variants combine shadows from the first batch with those from the
second batch through suffix channels in order to analyze the
overall set of game response strings for which any $y\in[n]$
survives. Theorem \ref{thm:PaulPackVol} gives winning conditions
for Paul in both variants in terms of a packing within a covering
of collections of these sets.  The nonexistence of such a packing
(resp., covering) provides a winning condition for Carole in the
original (resp., pathological) variant, in Theorem
\ref{thm:CarolePackVol} of Section \ref{sec:CaroleWinCond}.
Section \ref{sec:PaulCaroleProof} converts these winning
conditions into the main asymptotic results, Theorems
\ref{thm:PaulWinCond}-\ref{thm:CaroleWinCond}, after developing a
generalized Varshamov bound in Section \ref{sec:Varshamov}.
\begin{defn}\label{def:MrBalance}
Let $t\geq 2$ and $Q, M>0$ be integers, and let $r>0$.  A string
$w\in T^Q$ is {\em $(M,r)$-balanced} if, after splitting $w$ into
$M$ contiguous substrings of nearly equal length with (for
definiteness) length $\lceil Q/M\rceil$ substrings preceding
length $\lfloor Q/M\rfloor$ substrings, each letter of $T$ appears
in each section at most $\frac{1}{t}\lceil Q/M\rceil+r$ times. If
$w$ is not $(M,r)$-balanced, it is {\em $(M,r)$-unbalanced}.
Define $T^Q(M,r)$ to be the set of $(M,r)$-balanced strings of
$T^Q$.
\end{defn}

The Chernoff bound applies to the number of $(M,r)$-balanced
strings in $T^Q$.
\begin{lem}\label{lem:MrBalance}
Let $t\geq 2$ and
$Q,M>0$ be integers, and let
\begin{equation}\label{eqn:ri}
r(Q,M,i) =
\sqrt{\left\lceil\frac{Q}{M}\right\rceil\frac{\ln{(Mt2^i)}}{2}}.
\end{equation}
Then for $i\geq 1$, fewer than $t^Q2^{-i}$ strings in $T^Q$ are
$(M,r(Q,M,i))$-unbalanced.
\end{lem}
\begin{pf}
Select $w\in T^Q$ uniformly at random.  A fixed letter $a\in T$
appears independently with probability $1/t$ in each of the at
most $\lceil Q/M\rceil$ positions of a fixed section of $w$.
Letting $Y$ be the total number of occurrences of $a$ in this
section, by the standard Chernoff bound in Theorem A.1.4 of
\cite{AS00}, $\Pr(Y>\frac{1}{t}\lceil Q/M\rceil +r(Q,M,i)) <
\exp(-2r^2(Q,M,i)/\lceil Q/M\rceil) = \frac{2^{-i}}{Mt}$. The
result follows by subadditivity over $t$ letters and $M$
sections.\qed
\end{pf}

\subsection{Coding theoretic definitions\label{sec:codingTheoretic}}

Let $j,Q\geq 0$ and let $w=w_1\cdots w_Q,w'=w'_1\cdots w'_Q\in
T^Q$. The {\em Hamming ball} of radius $j$, or {\em $j$-ball},
centered at $w$ is the set $B(w,j) := \{w':0\leq d(w,w')\leq
j\}$, where $d(w,w'):=|\{i:w_i\neq w'_i\}|$ is the (Hamming)
distance between $w$ and $w'$.
We define distance to a set as usual; for example, for any
$w\in T^Q$, $d(w,T^Q(M,r)):=\min\{d(w,w'):w'\in T^Q(M,r)\}$ is
the distance between $w'$ and the set of $(M,r)$-balanced
strings of $T^Q$.
Our channel $C$ requires a generalization of Hamming balls to
$C$-shadows. Just as a $j$-ball is obtained from the center by
changing up to $j$ digits, a $C$-shadow is obtained from the
stem $w$ by applying a lie string $u\in C$ to $w$, as in
\eqref{eqn:applyLieString} for the case in which Paul's
questions are fully nonadaptive.
\begin{defn}
Let $w,w'\in T^Q$, and let $u\in\mathcal{L}(t)^j$. We write
$w\stackrel{u}{\rightarrow}w'$ if the lie string
$u=(a_1,b_1)\cdots(a_j,b_j)$ can be applied to $w$ to obtain $w'$
as in \eqref{eqn:applyLieString}.
\end{defn}

\begin{defn} \label{def:cShadow}
Let $w \in T^Q$, and let $C$ be a channel.  Then the {\em
$C$-shadow} $B(w,C)$ with stem $w$ is defined as
\begin{equation}
B(w,C) := \bigcup_{u\in C}\{w' \in T^Q :
w\stackrel{u}{\rightarrow}w'\}. \nonumber
\end{equation}
\end{defn}
Note that the $j$-ball $B(w,j)$ is a $C$-shadow with stem $w$ and
$C= \cup_{\ell=0}^j \mathcal{L}(t)^{\ell}$. A set
$\{B(w,C)\}_{(w,C)}$ of shadows is a {\em packing} in $T^Q$ if
$B(w,C)\cap B(w',C')=\emptyset$ for all distinct pairs of shadows in
the set.
The set is a {\em covering} of $T^Q$ if $\cup_{(w,C)}B(w,C)=T^Q$. An
$(x_i: 0 \leq i \leq k)$-packing (-covering) is a simultaneous
packing (covering) of $T^Q$ with $x_i$ $i$-balls for $0 \leq i \leq
k$. Similarly, an $(x_{C'}: C' \in \mathcal{I})$-packing (-covering)
is a simultaneous packing (covering) of $T^Q$ with $x_{C'}$
$C'$-shadows for $C'$ in some indexing set $\mathcal{I}$ of
channels. For our purposes, $\mathcal{I}=\mathcal{S}(C)$. From
coding theory, $A_t(Q,2j+1)$ is the maximum number of $j$-balls in a
packing of $T^Q$ (see, for example, \cite{PHB98}).
We define $b(Q,t,j)$ to be the size of any $j$-ball $B(w,j)$ in
$T^Q$, which is independent of $w$.  In particular,
\begin{equation}
    b(Q,t,j) = \sum_{\ell=0}^j\binom Q{\ell} (t-1)^{\ell}.
    \label{eqn:HammingBallSize}
\end{equation}
We make the following abbreviations in controlling $|B(w,C)|$ for
balanced $w$.
\begin{eqnarray}
G(Q,M,r,j) & := & \binom{M+j-1}{j}
       \left(\frac{1}{t}\left\lceil\frac{Q}{M}\right\rceil +
    r +k\right)^j, \qquad \mbox{and} \nonumber \\
H(Q,M,r,j) & := & \binom{M}{j}
       \left(
       \min\left(0,\frac{1}{t}\left \lceil \frac{Q}{M}
       \right \rceil-(t-1)r-2 -k \right) \right)^j.
\label{eqn:GH}
\end{eqnarray}
Here, $r$ corresponds to the balance tolerance parameter
$r(Q,M,i)$ in \eqref{eqn:ri}, and $k$ must appear to handle
stems $w$ with $d(w,T^Q(M,r))\leq k$.
\begin{lem}\label{lem:shadowBound}
Let $C$ be a channel,
let $w \in T^Q$ satisfy $d(w,s)\leq k$ for some
$(M,r)$-balanced $s\in T^Q$,
and let $u \in \mathcal{L}(t)^j$ be a lie string of length $j$.
Then
\begin{align}
H(Q,M,r,j) & \leq  |\{w':w \stackrel{u}{\rightarrow}{w'}\}|,
    |\{w':w'\stackrel{u}{\rightarrow}{w}\}|
    \leq  G(Q,M,r,j) \mbox{, and}
    \label{eqn:balWptoWBound}\\
\sum_{j=0}^{o(C)}\mathop{\sum_{u\in C,}}_{|u|=j} & H(Q,M,r,j) ~
\leq ~ |B(w,C)|
    ~ \leq ~ \sum_{j=0}^{o(C)}\mathop{\sum_{u\in C,}}_{|u|=j}
    G(Q,M,r,j)\,. \label{eqn:balShadBound}
\end{align}
\end{lem}
\begin{pf}
Divide $w$ and $s$ into $M$ contiguous sections as in
Defn.~\ref{def:MrBalance}. By definition, the maximum letter
frequency per section of an $(M,r)$-balanced $s\in T^Q$ is at
most $\frac{1}{t}\lceil Q/M\rceil +r$, and by subtraction the
minimum letter frequency per section is at least
$\frac{1}{t}\lceil Q/M\rceil -(t-1)r-2$.  Since $d(w,s)\leq k$,
add and subtract $k$ respectively from these quantities to get
corresponding bounds for $w$.
For \eqref{eqn:balWptoWBound}, we prove the bound on
$|\{w':w\stackrel{u}{\rightarrow}{w'}\}|$; the bound on
$|\{w':w'\stackrel{u}{\rightarrow}{w}\}|$ follows by replacing
$u=(a_1,b_1)\ldots(a_j,b_j)$ with $u'=(b_1,a_1)\ldots
(b_j,a_j)$ and noting that
$|\{w':w'\stackrel{u}{\rightarrow}{w}\}|=|\{w':w\stackrel{u'}
{\rightarrow}{w'}\}|$.  For the upper bound,
select
$j$ sections with possible repeats in $\binom{M+j-1}{j}$ ways,
to place the $j$ lies of $u$ in order. There are at most
$\left( \frac{1}{t}\lceil Q/M\rceil +r+k\right)$ ways of
applying lie $(a_{\ell},b_{\ell})$ within its section, so that
$w$ yields $w'$ as in \eqref{eqn:applyLieString}.
For the lower bound, use the bound $\frac{1}{t}\left \lceil Q/M
\right \rceil-(t-1)r-2 -k$ on the minimum letter frequency for
$w$, or 0 if this quantity is negative, and then
under-count by applying at most one lie per section.  A
$C$-shadow with stem $w$ consists of all $w'$ such that
$w\stackrel{u}{\rightarrow}w'$ for some $u\in C$. Thus
\eqref{eqn:balShadBound} follows by summing over $u$, graded by
length $|u|=j$, and applying \eqref{eqn:balWptoWBound}. \hfill
\qed
\end{pf}

We need the following lemma to handle applying lie strings to
severely unbalanced vertices in the pathological (original)
variant for Theorem \ref{thm:PaulPackVol}
(\ref{thm:CarolePackVol}). In fact, this is the motivation for
defining non-degenerate channels.

\begin{lem}
\label{lem:nondegen} Let $Q \geq t(k-1)+1$, and let $w \in T^{Q}$.
Let $C$ be a channel of order $k$.  If $C$ is non-degenerate with
respect to the original (pathological) variant, then $\bigcup_{u \in
C}  \{w': w \stackrel{u} {\rightarrow}w'\}$ $\left (\bigcup_{u \in
C}  \{w': w' \stackrel{u} {\rightarrow}w\}\right)$ is non-empty.
\end{lem}
\begin{pf}
Since $Q \geq t(k-1)+1$,
there exists a letter $c$ with minimum frequency $k$ in $w$.  In
the original (pathological) variant, let $a=c$ ($b=c$). By
Defn.~\ref{def:nondegen} there exists a $u\in C$ where $u=(a,b_1)
\cdots (a,b_j)$ ($u=(a_1,b) \cdots (a_j,b)$) for $0 \leq j \leq
k$. Construct $w'$ with $ w \stackrel{u} {\rightarrow}w'$ ($ w'
\stackrel{u} {\rightarrow}w$) by applying $u$ in $j$ arbitrarily
chosen positions in which $w$ has a $c$. \hfill \qed
\end{pf}

\subsection{A packing within covering condition for
Paul to win\label{sec:PaulWinCond}}

We now characterize a winning condition in
Thm.~\ref{thm:PaulWinCondPackCov} for Paul at the transition to
the second batch of questions, and go on to prove conditions
for which Paul can win the whole game in both variants in
Thm.~\ref{thm:PaulPackVol}. Paul's overall strategy in the
original variant is to split $[n]$ into blocks of size $\alpha$
and assign a unique block address, chosen from balanced
vertices in $T^{q_1}$, to each block. He forms his first batch
of $q_1$ questions by inspecting each element's block address.
Carole's first batch response $w'$ yields a state vector
$(x_{C'}(w'):C'\in\mathcal{S}(C))$, following
Defn.~\ref{def:stateVector}. Paul's selection of balanced block
addresses allows control on the entries of this state vector,
i.e., which $y\in [n]$ survive and in what fashion, through
Lemma \ref{lem:shadowBound}. He then wins the second batch of
$q_2$ questions, and thus the game (through the winning
strategy/packing equivalence in
Thm.~\ref{thm:PaulWinCondPackCov}), as follows.  He constructs
a packing in $T^{q_2}$ of the $C'$-shadows corresponding to
this state vector, fitting all but the singleton
$\{\varepsilon\}$-shadows inside Hamming balls centered on
balanced vertices of $T^{q_2}$ in order to ensure separation
and to control volume. The remaining empty space in $T^{q_2}$
exceeds the number of these singletons, and so Paul can add
singletons while preserving a packing, and his strategy is
winning. Paul's strategy in the pathological variant
piggy-backs his original variant strategy; he adds $\alpha'$
new elements to the above blocks of size $\alpha$ (thereby
increasing $n$). This increases the entries of
$(x_{C'}(w'):C'\in\mathcal{S}(C))$ enough so that the original
packing in $T^{q_2}$ can be augmented by new singletons to form
a covering. Unlike in the original variant, Paul must also
handle the case in which Carole's first batch response $w'$ is
not close to being balanced.  By assigning $t^{q_2}$ new
elements to each unbalanced block address in $T^{q_1}$, Paul
guarantees by virtue of non-degeneracy of the channel having
$t^{q_2}$ singletons to cover $T_{q_2}$. In either case he wins
by the winning strategy/covering equivalence in
Thm.~\ref{thm:PaulWinCondPackCov}.

\begin{thm} \label{thm:PaulWinCondPackCov}
Let $C$ be a channel, and let $(x_{C'}:C'\in\mathcal{S}(C))$ be
the state vector at the beginning of a 1-batch, $Q$-round game on
search space $[n]$. Then Paul
can win
the original (pathological) variant iff there exists an
$(x_{C'}:C'\in\mathcal{S}(C))$-packing (-covering) of $T^Q$.
\end{thm}
\begin{pf}
Given Paul's strategy in either variant and an element $y \in [n]$
counted by $x_{C'(y)}$, let $w(y) = w_1 \cdots w_Q$ be the
truthful response for $y$. For all $i$, $y$ is in the $w_i$th part
$A_{w_i}$ of Paul's $i$th question. Then $y$ survives the game iff
Carole responds with $w' \in B(w(y),C'(y))$. Paul wins the
original variant iff, for all responses $w'$, at most one $y$
survives, which occurs iff $\{B(w(y),C'(y))\}_{y \in [n]}$ is an
$(x_{C'}:C'\in\mathcal{S}(C))$-packing in $T^Q$. Similarly, Paul
wins the pathological variant iff for all responses $w'$, at least
one $y$ survives, which occurs iff $\{B(w(y),C'(y))\}_{y \in [n]}$
is an $(x_{C'}:C'\in\mathcal{S}(C))$-covering of $T^Q$. \hfill
\qed
\end{pf}

Empty sets are allowed in either the packing or covering of
Theorem \ref{thm:PaulWinCondPackCov}. Paul's strategy
determines the sets $B(w(y),C'(y))$, which might be empty when
$C'(y)$ violates Defn.~\ref{def:nondegen} or $Q$ is close to 0.
Adding empty sets neither hurts a packing nor helps form a
covering. In the next theorem, the parameters $M_1$ and $M_2$
are the number of sections into which the first and second
batches of $q_1$ and $q_2$ questions, respectively, are
divided, according to Defn.~\ref{def:MrBalance}. This
sectioning allows better counting of elements of $[n]$ with
particular game lie strings by considering the sections in
which lies occur.  The parameters $r_1$ and $r_2$ provide an
upper bound to the maximum letter frequency within sections in
the first and second batches of questions, respectively; and
$\eta_1$ and $\eta_2$ allow fine-tuning of $r_1$ and $r_2$ so
that an appropriately large proportion of the strings of
$T^{q_1}$ and $T^{q_2}$, respectively, are balanced.
Now we give the main conditions under which Paul has
winning strategies.
\begin{thm}\label{thm:PaulPackVol}
 Let $C$ be a $t$-ary channel of
order $k$, $q=q_1+q_2$ be the number of rounds split into two
positive integer batches, $\alpha,\alpha', M_1, M_2$ be
positive integers, and let $\eta_1,\eta_2$ be positive reals.
Following \eqref{eqn:ri}, define $r_1 := r(q_1,M_1, \eta_1
\log_2 q)$ and $r_2 := r(q_2, M_2, \eta_2 \log_2 q)$. Let
$c_k:=(k^2+3k-2)/2$, and define $G$ and $H$ as in
\eqref{eqn:GH}. If the packing condition
\begin{eqnarray}
\alpha\sum_{j=0}^i p_{k-i}^{(j)} G(q_1,M_1,r_1,j) & \leq &
    A_t(q_2-c_k,2(k-i)+1) -
     \frac{q^{-\eta_2} t^{q_2}}{b(q_2,t,k-i)}
    \label{eqn:packCond}
\end{eqnarray}
holds for all $1\leq k-i\leq k$, and the volume condition
\begin{eqnarray}
\alpha \sum_{i=0}^k\sum_{j=0}^i E_i(C)
 G(q_1,M_1,r_1,j)
    G(q_2,M_2,r_2,i-j)
    & \leq & t^{q_2}
    \label{eqn:volCond}
\end{eqnarray}
holds, then Paul can win the $(n,q_1,q_2,C)$-game when $n \leq
\alpha (1-q^{-\eta_1})t^{q_1}$.  Furthermore, if condition
\eqref{eqn:packCond} holds, $C$ is non-degenerate with respect to
the pathological variant, $q_1\geq t(k-1)+1$, and in addition the
volume condition
\begin{align}
\alpha \sum_{i=0}^k\sum_{j=0}^i &
    E_i(C) H(q_1,M_1,r_1,j)
    H(q_2,M_2,r_2,i-j)     \nonumber \\
    & + \alpha' \sum_{j=0}^k p_0^{(j)}  H(q_1,M_1,r_1,j)
    \ \geq \ t^{q_2}
    \label{eqn:pathVolCond}
\end{align}
holds, then Paul can win the $(n,q_1,q_2,C)^*$-game when $n\geq
(\alpha+\alpha')(1-q^{-\eta_1})t^{q_1}+t^{q}q^{-\eta_1}$.
\end{thm}

\begin{pf}
For the original variant,
Paul splits $[n]$ into blocks of size $\alpha$ and identifies
in bijective correspondence each block with an
$(M_1,r_1)$-balanced vertex of $T^{q_1}$.
Paul's first $q_1$ questions ask for the $q_1$ digits of the
distinguished element in this identification. By Lemma
\ref{lem:middleWinCond} and Theorem
\ref{thm:PaulWinCondPackCov}, Paul wins iff for any possible
answer $w'\in T^{q_1}$ by Carole, yielding the game state
vector $(x_{C'}(w'):C'\in\mathcal{S}(C))$ of
Defn.~\ref{def:stateVector}, there exists a
$(x_{C'}(w'):C'\in\mathcal{S}(C))$-packing in $T^{q_2}$.

{\it Claim 1.} \ If the packing condition \eqref{eqn:packCond}
holds, then there exists an $(x_{C'}(w'):C'\in\mathcal{S}(C),
o(C')\geq 1)$-packing in $T^{q_2}$ with all
stems in the set $\{z\in T^{q_2}:d(z,T^{q_2}(M_2,r_2))\leq
k\}$.
\smallskip

{\it Proof of Claim 1.} \ Because every order $(k-i)$ shadow fits
completely within a $(k-i)$-ball with the same center, it suffices
to show there exists an $(x_i:1\leq k-i\leq k)$-packing in $T^{q_2}$
with all centers
in $\{z\in T^{q_2}:d(z,T^{q_2}(M_2,r_2))\leq k\}$.
Let $\mathcal{D}_{k-i}\subseteq T^{q_2-c_k}$ be the set of
centers of a size $A_t(q_2-c_k,2(k-i)+1)$ packing of
$(k-i)$-balls.  We construct an
$(x_{k-i}=A_t(q_2-c_k,2(k-i)+1): 1 \leq k-i \leq k)$-packing of
Hamming balls in $T^{q_2}$ as follows:
$$
\begin{array}{clcccccc}
\mathcal{D} &= \quad &  \underbrace{1\cdots 1}_{k \ \mathrm{times}}
& 0\cdots 0 & 0\cdots 0 &\cdots & 00
    & \mathcal{D}_k\\
& \cup & 0\cdots 0 & \underbrace{1\cdots 1}_{k \ \mathrm{times}} &
0\cdots 0
    &\cdots & 00 & \mathcal{D}_{k-1}\\
& \cup & 0\cdots 0 & 0\cdots 0 & \underbrace{1\cdots 1}_{k-1 \
\mathrm{times}}
    &\cdots & 00 & \mathcal{D}_{k-2}\\
& & \vdots & \vdots & \vdots &\vdots &\vdots &\vdots \\
& \cup & 0\cdots 0 & 0\cdots 0 & 0\cdots 0
    &\cdots & \underbrace{11}_{2 \ \mathrm{times}} &
    \mathcal{D}_{1},
\end{array}
$$
in which the centers of the $(k-i)$-balls in $\mathcal{D}$ are taken
to be the extensions of their original centers in
$\mathcal{D}_{k-i}$. By construction, two distinct balls of radius
$i$ and $j$ are disjoint, as the distance between their centers is
at least $i+j+1$.
By Lemma \ref{lem:MrBalance}, for fixed $k-i$, the number of
$(k-i)$-balls in the packing $\mathcal{D}$ comprised entirely
of $(M_2,r_2)$-unbalanced vertices is at most
$q^{-\eta_2}t^{q_2}/b(q_2,t,k-i)$. After deleting the
corresponding centers from $\mathcal{D}$, at least
$A_t(q_2-c_k,2(k-i)+1)-q^{-\eta_2}t^{q_2}/b(q_2,t,k-i)$
$(k-i)$-balls with centers
in $\{z\in T^{q_2}:d(z,T^{q_2}(M_2,r_2))\leq k\}$
remain. We finish the claim by showing that for any response
$w'$ by Carole, the left-hand side of \eqref{eqn:packCond} is
an upper bound on $x_i$ for $1\leq k-i\leq k$:
\begin{eqnarray}
x_i & =&\mathop{\sum_{C'\in\mathcal{S}(C)}}_{o(C')=k-i}x_{C'} \leq
    \alpha\mathop{\sum_{C'\in\mathcal{S}(C)}}_{o(C')=k-i}
    \sum_{u,S_C(u)=C'} |\{w:w\stackrel{u}
    {\rightarrow}w'\}| \label{eqn:xiToxc} \\
& = & \alpha \sum_{j=0}^i\mathop{\sum_{u,|u|=j}}_{o(S_C(u))=k-i}
    |\{w:w \stackrel{u}{\rightarrow}w'\}|
    \label{eqn:xcexact}
 \ \leq \ \alpha\sum_{j=0}^i p_{k-i}^{(j)} G(q_1,M_1,r_1,j)
    . \label{eqn:applUBXi}
\end{eqnarray}
Line \eqref{eqn:xiToxc} is by the definitions of $x_i$, $x_{C'}$,
and $w\stackrel{u}{\rightarrow}w'$, where the inequality is because
the unbalanced $w$ contribute nothing to $x_{C'}$.  The equality in
\eqref{eqn:xcexact} is by a straightforward reindexing of the
summation.  The inequality in \eqref{eqn:xcexact}
is by the definition of $p^{(j)}_{k-i}$ and Lemma
\ref{lem:shadowBound} since
Carole's response $w'$ must satisfy $d(w',T^{q_1}(M_1,r_1))\leq
k$, or else she defaults as
all elements of $[n]$ are identified with $(M_1,r_1)$-balanced
strings.

{\it Claim 2.} \ If the packing condition \eqref{eqn:packCond} and
the volume condition \eqref{eqn:volCond} both hold, then there
exists an $(x_{C'}(w'):C'\in\mathcal{S}(C))$-packing in
$T^{q_2}$.\smallskip

{\it Proof of Claim 2.}\ To show that such a packing exists, it
is enough to demonstrate that $x_k$ 0-balls (singletons) can be
packed in the unoccupied space of the packing in Claim 1.  We
therefore show that the total volume of the Claim 1 packing and
the singletons is at most the left hand side of
\eqref{eqn:volCond}. For $0 \leq i <k$, every $S_C(u)$-shadow
in the packing in $T^{q_2}$ counted by \eqref{eqn:xcexact} is
bounded in size by Lemma \ref{lem:shadowBound}, because all
stems of the packing in Claim 1 are in $\{z\in
T^{q_2}:d(z,T^{q_2}(M_2,r_2))\leq k\}$.
Hence the space occupied by the packing from Claim 1 is
\begin{equation} \leq
\alpha \sum_{i=0}^k \sum_{j=0}^i
    \mathop{\sum_{u,|u|=j}}_{o(S_C(u))=k- i} |\{w:w\stackrel{u}
    {\rightarrow}w'\}|
    \mathop{\max_{z\in T^{q_2}}}_{d(z,T^{q_2}(M_2,r_2))\leq k}
    |B(z,S_C(u))|.
    \label{eqn:rawVolBound}
\end{equation}
The singletons, counted when $i=k$, are all volume 1
$\{\epsilon\}$-shadows with $B(z,\{\epsilon \}) = \{z\}$
regardless of $z$.  Applying Lemma \ref{lem:shadowBound},
\eqref{eqn:rawVolBound} is
\begin{eqnarray}
& \leq & \alpha \sum_{i=0}^k \sum_{j=0}^i
\mathop{\sum_{u,|u|=j}}_{o(S_C(u)) =k- i} G(q_1,M_1,r_1,j)
   \sum_{\ell=0}^{k-i} \,\mathop{\sum_{v \in
    S_C(u)}}_{|v|=\ell}
   G(q_2,M_2,r_2,\ell)  \label{eqn:mess1}\\
& = &\alpha \sum_{i=0}^k \sum_{j+\ell \leq k}\sum_{\substack{uv \in
    E_{j+\ell}(C)  \\ o(S_C(u))=k-i \\ |u|=j, |v|=\ell}}
     G(q_1,M_1,r_1,j)   G(q_2,M_2,r_2,\ell) \label{eqn:mess2}\\
& =&\alpha\sum_{j+\ell \leq k}\,\sum_{\substack{uv \in E_{j+\ell}(C)    \\
|u|=j, |v|=\ell}}
     G(q_1,M_1,r_1,j)    G(q_2,M_2,r_2,\ell) \label{eqn:mess3}\\
& = & \alpha\sum_{i=0}^k \sum_{j=0}^i E_i(C)
     G(q_1,M_1,r_1,j)    G(q_2,M_2,r_2,i-j). \label{eqn:mess4}
 \end{eqnarray}

The double sums over $j$ and $\ell$ combine to a sum over $j+\ell
\leq k$ in \eqref{eqn:mess2} because, for $j>i$, there are no $u$
with $|u|=j$ and $o(S_C(u))=(k-i)$.  By interchanging the first two
summations in \eqref{eqn:mess2}, the sum over $i$ has the effect of
summing over the orders of $S_C(u)$, and since each $u \in C$ has a
unique $o(S_C(u))$, we have \eqref{eqn:mess3}. Finally, setting
$j+\ell = i$ gives \eqref{eqn:mess4}, completing the proof of Claim
2.

Therefore Paul wins if $n$ is $\alpha$ times the number of
$(M_1,r_1)$-balanced vertices in $T^{q_1}$, which is at least
$(1-q^{-\eta_1})t^{q_1}$ by Lemma \ref{lem:MrBalance}. If $n$ is
less than this number, Paul can clearly still win by simply removing
shadows from the packing.

For the pathological variant, Paul identifies $(\alpha +
\alpha')$ elements of $[n]$ with each of the
$(M_1,r_1)$-balanced vertices of $T^{q_1}$, and $t^{q_2}$
elements
with
each $(M_1,r_1)$-unbalanced vertex. We may assume $\alpha \leq
\alpha + \alpha' <t^{q_2}$, for suppose to the contrary $\alpha
+ \alpha' \geq t^{q_2}$.
Then Paul can win by the following argument.
Let $w'$ be Carole's response after the first batch of
questions. By Lemma \ref{lem:nondegen}, there exist  $u \in C$
and $w \in T^{q_1}$ with $\{w \stackrel{u} {\rightarrow}w'\}$.
There are at least $t^{q_2}$ elements identified with $w$
that will survive the first
batch with suffix channel $C' = S_C(u)$ containing $\epsilon$.
In the second batch, Paul identifies at least one of these
elements counted by $x_{C'}$ to each vertex of $T^{q_2}$ and
asks for the $q_2$ digits of the distinguished element.
Regardless of Carole's response $z' \in T^{q_2}$, the
element(s) identified with $z'$ survives.

As in the original variant, Paul's first $q_1$ questions ask for the
$q_1$ digits of the distinguished element in the above
identification.  By Theorem \ref{thm:PaulWinCondPackCov}, Paul wins
iff for all possible answers $w'\in T^{q_1}$ by Carole, yielding the
game state vector $(x_{C'}(w'):C'\in\mathcal{S}(C))$ of
Defn.~\ref{def:stateVector}, there exists an
$(x_{C'}(w'):C'\in\mathcal{S}(C))$-covering of $T^{q_2}$.

{\it Claim 3.} \  If the packing condition \eqref{eqn:packCond}
holds, $C$ is non-degenerate with respect to the pathological
variant, $q_1\geq t(k-1)+1$, and the volume condition
\eqref{eqn:pathVolCond} holds, then there exists an
$(x_{C'}(w'):C'\in\mathcal{S}(C))$-covering of $T^{q_2}$.

{\it Proof of Claim 3.} \ We first consider the case
in which $d(w',T^{q_1}(M_1,r_1))\leq k$.
For this case, consider only $\alpha+\alpha'$ elements of $[n]$
identified per $w\in T^{q_1}$.  The state vector
$(x_{C'}(w'):C'\in\mathcal{S}(C))$ can be split into
$(x_{C'}(w')|_{\alpha}:C'\in\mathcal{S}(C))$ and
$(x_{C'}(w')|_{\alpha'}:C'\in\mathcal{S}(C))$ from the
contributions of $\alpha$ and $\alpha'$ respectively. By
\eqref{eqn:packCond} and Claim 1, there exists an
$(x_{C'}(w')|_{\alpha}:C'\in\mathcal{S}(C), o(C')\geq
1)$-packing in
$T^{q_2}$, with all stems in the set $\{z\in
T^{q_2}:d(z,T^{q_2}(M_2,r_2))\leq k\}$.
If the volume after adding in the $\alpha$ elements identified
with each singleton is $>t^{q_2}$ we have a covering, and we
are done.  Otherwise, we have an
$(x_{C'}(w')|_{\alpha}:C'\in\mathcal{S}(C))$-packing in
$T^{q_2}$.

Considering only the first $\alpha$ elements identified with each
vertex of $T^{q_1}$, equation \eqref {eqn:xiToxc} becomes equality
throughout, and the volume of the
$(x_{C'}(w')|_{\alpha}:C'\in\mathcal{S}(C))$-packing (analogous to
\eqref{eqn:rawVolBound}) is
\begin{equation}
\geq \sum_{i=0}^k \sum_{j=0}^i
    \mathop{\sum_{u,|u|=j}}_{o(S_C(u)) =k- i}\alpha
    |\{w:w\stackrel{u} {\rightarrow}w'\}|
    \mathop{\min_{z\in T^{q_2}}}_{d(z,T^{q_2}(M_2,r_2))\leq k}
    |B(z,S_C(u))|.
    \label{eqn:lowerRawVolBound}
\end{equation}
The
singletons, counted when $i=k$, are all volume 1
$\{\epsilon\}$-shadows with $B(z,\{\epsilon \}) = \{z\}$
regardless of $z$.
 Applying Lemma
\ref{lem:shadowBound}, and manipulating the summations as in
\eqref{eqn:mess1}--\eqref{eqn:mess4}, \eqref{eqn:lowerRawVolBound}
is
\begin{eqnarray}
& \geq & \alpha \sum_{i=0}^k \sum_{j=0}^i E_i(C)
    H(q_1,M_1, r_1, j)H(q_2,M_2,r_2,i-j).
    \label{eqn:initialCover}
\end{eqnarray}
Following \eqref{eqn:xiToxc}, the additional $\alpha'$ elements
identified with each vertex of $T^{q_1}$ have a contribution to
$x_{\{\epsilon\}}$ of exactly
\begin{eqnarray}
x_{\{\epsilon\}}|_{\alpha'}& = &\sum_{u,S_C(u)=\{\epsilon\}}
    \alpha' |\{w: w
    \stackrel{u}{\rightarrow}w'\}| \ \geq \ \sum_{j=0}^k
    \mathop{\sum_{u,|u|=i}}_{S_C(u)=\{\epsilon\}} \alpha'
    H(q_1,M_1,r_1,j) \nonumber \\
& =& \alpha' \sum_{j=0}^k p^{(j)}_0 H(q_1, M_1, r_1, j) \nonumber
\end{eqnarray}
which by volume condition \eqref{eqn:pathVolCond} and
\eqref{eqn:initialCover} is at least the number of vertices of
$T^{q_2}$ not covered in the
$(x_{C'}(w')|_{\alpha}:C'\in\mathcal{S}(C))$-packing. We extend to
an $(x_{C'}(w'):C'\in\mathcal{S}(C))$-covering of $T^{q_2}$ by using
at most $x_{\{\epsilon\}}|_{\alpha'}$ $\{\epsilon\}$-shadows to
cover, with $B(z,\{\epsilon\})$, any $z \in T^{q_2}$ not covered by
the packing. Based on the covering constraint, any unaccounted-for
shadows in the original identification of elements of $[n]$ to
vertices of $T^{q_1}$ may be ignored.

Now assume that Carole responds to the first $q_1$ questions with
$w'\in T^{q_1}$
and $d(w',T^{q_1}(M_1,r_1))>k$.
Again, let $(x_{C'}(w'):C'\in\mathcal{S}(C))$ be the state
vector after Carole's response.  By Lemma \ref{lem:nondegen},
there exists a $u\in C$ and a $w\in T^{q_1}$ such that
$w\stackrel{u}{\rightarrow}w'$. Since $|u|\leq k$, $d(w',w)\leq
k$, and so $w$ is $(M_1,r_1)$-unbalanced. Thus the $t^{q_2}$
elements identified with $w$ have suffix channel $S_C(u)$
containing $\epsilon$, so that $x_{S_C(u)}\geq t^{q_2}$. The
collection of $S_C(u)$-shadows $\{B(z,S_C(u)):z\in T^{q_2}\}$
covers $T^{q_2}$, since $\epsilon\in S_C(u)$ implies $z\in
B(z,S_C(u))$. There exists an
$(x_{C'}(w'):C'\in\mathcal{S}(C))$-covering of $T^{q_2}$ by
placing the remaining shadows arbitrarily.

Therefore, whether $w'$ is
close to being
balanced or not, Paul wins with $n$ equal to $(\alpha+\alpha')$
times the number of $(M_1,r_1)$-balanced vertices of $T^{q_1}$
and $t^{q_2}$ times the number of $(M_1,r_1)$-unbalanced
vertices of $T^{q_1}$.  By Lemma \ref{lem:MrBalance}, this is
at most
$(\alpha+\alpha')(1-q^{-\eta_1})t^{q_1}+t^{q}q^{-\eta_1}$. Paul
can win for any $n'>n$ by treating extra elements arbitrarily
without disturbing the covering constructed above. \hfill \qed
\end{pf}

\subsection{A condition for Carole to win\label{sec:CaroleWinCond}}

The following theorem gives conditions under which Carole has
winning strategies in both the original and pathological 2-batch
games, by way of nonexistence of packings or coverings,
respectively, corresponding to winning strategies for Paul.

\begin{thm}\label{thm:CarolePackVol}
Let $C$ be a $t$-ary channel of order $k$, $q=q_1+q_2$ be the
number of rounds split into two positive integer batches, $M_1,
M_2$ be positive integers, and let $\eta_1,\eta_2,\eta$ be
positive reals. Following \eqref{eqn:ri}, define $r_1 :=
r(q_1,M_1, \eta_1 \log_2 q)$, $r_2 := r(q_2, M_2, \eta_2 \log_2
q)$, and $r := r(q, M, \eta \log_2 q)$; define $G$ and $H$ as
in \eqref{eqn:GH}. If $C$ is non-degenerate with respect to the
original variant, $q_1,q_2\geq t(k-1)+1$, and, in addition, the
volume condition
\begin{eqnarray}
n & > & t^q\big( E_k(C)\left(H(q_1,M_1,r_1,k)+H(q_2,M_2,r_2,k)\right)
    \big)^{-1} \nonumber\\
& &  + (q^{-\eta_1}  + q^{-\eta_2})t^q \label{eqn:CaroleVolCond}
\end{eqnarray}
holds, then Carole can win the $(n,q_1,q_2,C)$-game.  If the volume
condition
\begin{eqnarray}
 n & < & t^q\left(1-q^{-\eta}\right) \left(\sum_{i=0}^k
    E_i(C)\binom{M+i-1}{i}\left(\frac1t\left\lceil
    \frac qM \right\rceil + r + 1\right)^i\right)^{-1}
\label{eqn:PathCaroleVolCond}
\end{eqnarray}
holds, then Carole can win the $(n,q_1,q_2,C)^*$-game.
\end{thm}

\begin{pf}
We define the {\em response set} of an element $y$ of the search
space $[n]$ to be $\mathcal{R}(y):=\{w'z' \in T^q: y\mbox{ survives
with game response string } w'z'\}$.  Call a response string $w'z'$
{\em doubly balanced} if $w' \in T^{q_1}$ is $(M_1,r_1)$-balanced,
and $z' \in T^{q_2}$ is $(M_2,r_2)$-balanced. We say $y$ is {\em
typical} if every $w'z'
\in \mathcal{R}(y)$ is doubly balanced.

For the original variant, Carole can win if $\{ \mathcal{R}(y): y
\in [n]\}$ is not a packing, since if $\mathcal{R}(y) \cap
\mathcal{R}(y') \neq \emptyset$, then there exists a $w'z'$ for
which both $y$ and $y'$ survive.  Assume $y$ is typical. Let
$w=w(y)$ be the truthful first batch response for $y$. We may
assume
$d(w,T^{q_1}(M_1,r_1))\leq k$.
Otherwise, since $C$ is non-degenerate, by Lemma
\ref{lem:nondegen} there exists a $u \in C$ such that $w
\stackrel{u}{\rightarrow}{w'}$ with $w'$
$(M_1,r_1)$-unbalanced, making $y$ atypical.

We under-count $\mathcal{R}(y)$ by counting only those $w'$
with $w \stackrel{u}{\rightarrow}{w'}$ for $u=\epsilon$ or
$u\in C$ with $|u|=k$, thus guaranteeing that $S_C(u)$ is
non-degenerate. Using Lemma \ref{lem:shadowBound}, the number
of such $w'$ is at least $1+E_k(C)H(q_1,M_1,r_1,k)$, where the
1 term corresponds to $u=\epsilon$.  If $y$ survives the first
batch given that Carole's response is one of these $w'$, then
Paul's strategy determines the truthful $z$ for $y$ in the
second batch.  As before, $z$ must
satisfy $d(z,T^{q_2}(M_2,r_2))\leq k$,
or else $y$ survives for an unbalanced $z'$.  The number of
$z'$ for which $y$ survives the second batch is dependent on
the suffix channel $S_C(u)$, and is at least
$E_k(C)H(q_2,M_2,r_2,k)$ when $u=\epsilon$ and 1 otherwise.
Therefore, similar to \eqref{eqn:mess1}--\eqref{eqn:mess4}
(with many terms omitted), the size of $\mathcal{R}(y)$ is at
least
\begin{equation}
E_k(C)\left( H(q_1,M_1,r_1,k) + H(q_2,M_2,r_2,k)\right).
\nonumber
\end{equation}

If $y$ is atypical, then there exists a response sequence $w'z'$ for
which $y$ survives, and either $w'$ is $(M_1,r_1)$-unbalanced, or
$z'$ is $(M_2,r_2)$-unbalanced.  Thus there is at least one
non-doubly balanced string in $\mathcal{R}(y)$.   By Lemma
\ref{lem:MrBalance}, there are at most $(q^{-\eta_1} +
q^{-\eta_2})t^q$ such $w'z'$.

We can pack at most $t^q \left(E_k(C)\left( H(q_1,M_1,r_1,k) + H(q_2,M_2,r_2,k)\right)\right)^{-1}$ response
    sets for $y$ typical and, independently, at most $(q^{-\eta_1} +
q^{-\eta_2})t^q$ response sets for $y$ atypical.  Therefore, if
\eqref{eqn:CaroleVolCond} holds, Carole can win.

For the pathological variant, Carole can win if $\{ \mathcal{R}(y):
y \in [n]\}$ is not a covering, since for $w\in T^q$, if $w' \notin
\cup_{y\in[n]}\mathcal{R}(y)$, then $w'$ is a response for which no
element $y$ survives. We further handicap Carole by allowing Paul
full adaptivity, i.e., Paul can wait to ask each question until
after Carole responds to the previous question.  If Carole can win
the fully adaptive case, she can certainly win the two-batch case
for any $q_1$ and $q_2$.

We bound the number of $(M,r)$-balanced strings in
$\mathcal{R}(y)$ for arbitrary $y \in [n]$. Each balanced string
in $\mathcal{R}(y)$ is a result of applying a length $i$ lie
string, $0 \leq i \leq k$, to the truthful sequence of responses
to Paul's queries, and is therefore identified by the lie string
and positions of the lies. Divide Carole's $q$ responses into $M$
blocks as in Defn.~\ref{def:MrBalance}.
Carole selects $u
 \in C$ with $|u|=i$ in
$E_i(C)$ ways, and the $i$ sections in which to place the lies of
$u$ in order in $\binom{M+i-1}{i}$ ways.

The first lie $(a,b)$ of $u$ to be placed in any block must occur
within the first $\left (\frac{1}{t} \lceil \frac{q}{M}\rceil + r
+1\right)$ occurrences of $a$ in that block; otherwise all of the
resulting game response strings will be unbalanced.
This restriction holds for every subsequent lie; therefore the
maximum number of balanced strings in $\mathcal{R}(y)$ is at most
\begin{equation}
\label{eqn:CaroleNonAdaptPath}
 \sum_{i=0}^{k} E_i(C) {\binom{M+i-1}{i} }\left (\frac{1}{t} \left \lceil
\frac{q}{M} \right \rceil + r +1 \right)^i. \nonumber
 \end{equation}
By Lemma \ref{lem:MrBalance}, there are at least $t^q\left(1
-q^{-\eta }\right)$ $(M,r)$-balanced strings in $t^q$.  Thus at
least $t^q\left(1 -q^{-\eta} \right) \left(\sum_{i=0}^{k} E_i(C)
{\binom{M+i-1}{i} }\left (\frac{1}{t} \left \lceil \frac{q}{M}
\right \rceil + r +1 \right)^i \right)^{-1}$ response
    sets are necessary to cover $t^q$.  Therefore, if
\eqref{eqn:PathCaroleVolCond} holds, Carole can win.
\hfill \qed
\end{pf}

\subsection{An asymptotic approximation and Varshamov bound for
the main theorem\label{sec:Varshamov}}

In order to convert Theorems
\ref{thm:PaulPackVol}-\ref{thm:CarolePackVol} to asymptotic
form, we require the technical Lemma \ref{lem:MReduct} and an
asymptotic form of a generalized Varshamov bound in Corollary
\ref{cor:asympGenVarshamov}. Lemma \ref{lem:MReduct} will be
used several times to approximate quantities such as $G$ and
$H$ of \eqref{eqn:GH}.  An asymptotic version of the packing
condition \eqref{eqn:packCond} is allowed by bounding
$A_t(Q,2R+1)$ with Theorem \ref{thm:Varshamov}, when $t$ is a
prime power, generalizing to $t$ not a prime power with Lemma
\ref{lem:alphabetChange}, and converting to asymptotic form
with Corollary \ref{cor:asympGenVarshamov}.

\begin{lem}\label{lem:MReduct}
Let $\ell\in\mathbb{Z}$, $j\in\mathbb{Z}_{\geq 0}$,
$\constG,\constH\in\mathbb{R}$, and $\eta\in\mathbb{R}^{+}$ be constants.
Let $q\rightarrow\infty$, and let $f(q)$ be nonnegative with
$f(q)\rightarrow\infty$. Let $Q$ satisfy $(\ln q)^{3/2}f(q)\leq
Q\leq q$. Let $M=\lceil Q^{1/3}\rceil$.  Then
\begin{align}
\binom {M+\ell}{j} & \left(\frac 1t \left\lceil\frac
    {Q}{M}\right\rceil +
    \constG\sqrt{\left\lceil\frac{Q}{M}\right\rceil
    \frac{\ln(tM q^{\eta})}{2}} +\constH\right)^j\nonumber\\
& = \binom{Q}{j}t^{-j}\left(1 +\frac{tj\constG}{\sqrt{2}}
    \frac{\sqrt{\ln(Q^{1/3}q^\eta)}}{Q^{1/3}}(1+o(1))\right)\,.
    \label{eqn:MReduct}
\end{align}
\end{lem}
\begin{pf}
First, note that for any $N\rightarrow\infty$,
\begin{equation}\label{eqn:binomReduct}
    \binom{N+\ell}{j} = \frac{N^j}{j!}\left(1\pm
    \Theta\left(\frac{1}{N}\right)\right).
\end{equation}
Applying \eqref{eqn:binomReduct} to $\binom{M+\ell}{j}$ and
replacing $\lceil Q/M\rceil$ with $Q/M+\Theta(1)$, the left-hand
side of \eqref{eqn:MReduct} becomes
\begin{align}
= & \frac{M^j}{j!}\left(1\pm\Theta\left(\frac{1}{M}\right)\right)
    \left(\frac{Q}{tM} + \constG\sqrt{\left(\frac{Q}{M}+O(1)\right)
    \frac{\ln(tMq^{\eta})}{2}} \pm O(1)\right)^j.\nonumber
\end{align}
Factoring out $(Q/(tM))^j$ from the last factor, and expanding the
square root and the entire last factor using $(1+x)^p=1+px+O(x^2)$
as $x\rightarrow 0$, this becomes
\begin{align}
= & \frac{M^j}{j!}\frac{Q^j}{t^jM^j}
    \left(1\pm\Theta\left(\frac{1}{M}\right)\right)
    \left(1 + tj\constG\sqrt{\frac{M}{Q}
    \frac{\ln(tMq^{\eta})}{2}} + O\left(\frac{M}{Q}\ln q
    \right)\right).\nonumber
\end{align}
Applying \eqref{eqn:binomReduct} to $Q^j/j!$, noting that
$\ln(t)=o(\ln(q))$ and $M=Q^{1/3}+O(1)$, and applying the binomial
expansion on the square root again, this becomes the right-hand
side of \eqref{eqn:MReduct}.  The constraint $Q\geq (\ln
q)^{3/2}f(q)$ allows collection of lower order terms into the
$(1+o(1))$ factor.
\hfill \qed
\end{pf}

The Varshamov lower bound for $A_t(Q,2R+1)$, for $t$ a prime power,
relies on a construction on a vector space over a finite field, and
can be found, for example, as Theorem 3.4 of \cite{PHB98}.  A {\em
linear code} with minimum distance $2R+1$ may be viewed as a packing
of $R$-balls whose centers form a vector space.
\begin{thm}[Varshamov bound]\label{thm:Varshamov} Let $t\geq 2$ be
a prime power, and let $Q\geq 1$ and $R\geq 0$ be integers.  Then
$$
A_t(Q,2R+1)\geq B_t(Q,2R+1) \ \geq \
    t^{Q-\left\lceil \log_t \left(1+\sum_{i=0}^{2R-1}
    \binom{Q-1}{i}(t-1)^i\right)\right\rceil}.
$$
where $B_t(Q,2R+1)$ is the maximum size of a linear code of length
$Q$, alphabet $t$, and minimum distance $2R+1$.
\end{thm}
The following lemma allows extension to a weakened version of the
Varshamov bound for $t$ not a prime power.  The lemma is due to
Gevorkyan \cite{G79}, but can also be found within the proof of
Lemma 2 of \cite{YD04}.  We include the proof here for clarity.
\begin{lem}\label{lem:alphabetChange}
Let $t_2\geq t_1\geq 2$, $Q\geq 1$ and $R\geq 0$ be integers. Then
$$
A_{t_1}(Q,2R+1) \ \geq \ (t_1/t_2)^Q A_{t_2}(Q,2R+1).
$$
\end{lem}
\begin{pf}
Define $T_1:=\{0,\ldots,t_1-1\}$ and $T_2:=\{0,\ldots,t_2-1\}$.  Let
$\mathcal{V}_2\subseteq T_2^Q$ be the set of strings which are the
centers of $R$-balls in a packing of $T_2^Q$. In particular, assume
$|\mathcal{V}_2|=A_{t_2}(Q,2R+1)$. View the strings of
$\mathcal{V}_2\subseteq  T_2^Q$ as elements of the additive group
$\mathbb{Z}_{t_2}^Q$.  Let $w\in  T_2^Q$, and define the translation
$\mathcal{V}_2^w:=\{w+w':w'\in \mathcal{V}_2\}$ of $\mathcal{V}_2$.
Let $z\in \mathbb{Z}_{t_1}^Q$ be viewed as an element of
$\mathbb{Z}_{t_2}^Q$ since $t_1\leq t_2$. Since $\mathbb{Z}_{t_2}^Q$
is a group, the number of translations of $\mathcal{V}_2$ containing
$z$ is $|\{w:z\in \mathcal{V}_2^w\}|=|\mathcal{V}_2|$. The average
number of elements of $\mathbb{Z}_{t_1}^Q$ contained in a translate
$\mathcal{V}_2^w$ is thus $|\mathbb{Z}_{t_1}^Q|\cdot
|\mathcal{V}_2|/|\mathbb{Z}_{t_2}^Q|=(t_1/t_2)^Q|\mathcal{V}_2|$,
and there must exist a translate $\mathcal{V}_2^{w^*}$ with at least
this average. Translation preserves distance between strings, and so
we may take the centers of the $R$-balls in our packing of $T_1^Q$
to be $\mathcal{V}_2^{w^*}\cap T_1^Q$.
\hfill \qed
\end{pf}

\begin{cor}\label{cor:asympGenVarshamov}
Let $t\geq 2$, $Q\geq 1$ and $k\geq R\geq 0$ be integers. Let
$c_k=(k^2+3k-2)/2$ as in Theorem \ref{thm:PaulPackVol}, and let
$\constI < 1$ be an arbitrary constant. Then with $t$ and $R$ fixed,
for $Q$ large enough,
$$
A_t(Q-c_k,2R+1) \ > \
    \constI \frac{(2R-1)!}{2^{2R}(t-1)^{2R-1}}
    \frac{t^{Q-c_k-1}}{Q^{2R-1}}.
$$
\end{cor}
\begin{pf}
Set $t_1=t$, and let $t_2$ be the smallest prime power for which
$t_1\leq t_2$.  In particular, $t_2\leq 2t_1-1$. Applying Lemma
\ref{lem:alphabetChange} and Theorem \ref{thm:Varshamov}, we have
\begin{align}
A_{t_1}&(Q-c_k,2R+1) \ \geq \ t_1^{Q-c_k} t_2^{-\left\lceil
\log_{t_2}
    \left(1+\sum_{i=0}^{2R-1}
    \binom{Q-c_k-1}{i}(t_2-1)^i\right)\right\rceil} \nonumber \\
& \geq \ \frac{t_1^{Q-c_k}t_2^{-1}}
    {1+\sum_{i=0}^{2R-1} \binom{Q-c_k-1}{i}(t_2-1)^i}
\ > \ \frac12 \frac{t^{Q-c_k-1}}
    {1+\sum_{i=0}^{2R-1} \binom{Q-c_k-1}{i}2^i(t-1)^i}\,, \nonumber
\end{align}
from which the result follows by observing that the denominator is
dominated by the $i=2R-1$ term.
\hfill \qed
\end{pf}

\subsection{Proofs of the main result: winning conditions
for Paul and Carole\label{sec:PaulCaroleProof}}

We now apply the results of Section \ref{sec:Varshamov} and
standard asymptotic analysis to prove Theorems
\ref{thm:PaulWinCond}-\ref{thm:CaroleWinCond}.

{\bfseries{Proof of Theorem \ref{thm:PaulWinCond}.}}
Let $\constB>\constJ > tk(k+1)\sqrt{k+2}/\sqrt{2}$ be constants and let
$$
\alpha = \left\lfloor \frac{t^{q_2+k}}{E_k(C)\binom{q}{k}}
    \left(1-\constJ\frac{\sqrt{\ln q}}{q_2^{1/3}}\right)\right\rfloor.
$$
Fix $M_1=\lceil q_1^{1/3}\rceil$ and $M_2=\lceil
q_2^{1/3}\rceil$. Let $\eta_1=\eta_2=k+1$. Let
$\constA<\constK<\left({\frac{(2k-1)!t^{-c_k-k-1}E_k(\mathcal{C})}
{2^{2k}(t-1)^{2k-1}k!}}\right)^{1/(2k-1)}$ be constants, and
let $\constL>1$ and $\constI<1$ be constants such that
$\constA<(\constI/\constL)^{1/(2k-1)}\constK$.  Applying Lemma
\ref{lem:MReduct} and \eqref{eqn:binomReduct}, and noting that
the terms for $0\leq j<i$ are asymptotically negligible, for
$q$ (and thus $q_1$) sufficiently large, the left-hand side of
\eqref{eqn:packCond} is
\begin{align}\label{eqn:lhsPackCondAsymp}
    < & \constL \frac{t^{q_2+k}k!}{E_k(C)q^k}
    \frac{p_{k-i}^{(i)}}{t^i}\frac{q_1^i}{i!}.
\end{align}
Applying Corollary \ref{cor:asympGenVarshamov} and
\eqref{eqn:HammingBallSize}, for $q$ (and thus $q_2$) sufficiently
large, the right-hand side of \eqref{eqn:packCond} is
\begin{equation}\label{eqn:rhsAsympPacking}
A_t(q_2-c_k,2(k-i)+1) -
     \frac{q^{-k-1} t^{q_2}}{b(q_2,t,k-i)}
    >  \ \constI \frac{(2(k-i)-1)!}{2^{2(k-i)}(t-1)^{2(k-i)-1}}
        \frac{t^{q_2-c_k-1}}{q_2^{2(k-i)-1}},
\end{equation}
since the assumption $1\leq k-i\leq k$ for \eqref{eqn:packCond}
makes $q^{-k-1} t^{q_2}/b(q_2,t,k-i)$ asymptotically
negligible. For \eqref{eqn:packCond} to hold, it will suffice
that \eqref{eqn:lhsPackCondAsymp} be at most the right-hand
side of \eqref{eqn:rhsAsympPacking}, for
$1\leq k-i\leq k$. This is immediate when $p^{(i)}_{k-i}= 0$;
otherwise the condition is equivalent to
\begin{equation}
q_2^{2(k-i)-1} \ \leq \ \frac{\constI}{\constL}\frac{(2(k-i)-1)!
    t^{-c_k-(k-i)-1}E_k(C)i!}{2^{2(k-i)}(t-1)^{2(k-i)-1}p_{k-i}^{(i)} k!}
    \frac{q^k}{q_1^i}
    . \label{eqn:alphaRootBound}
\end{equation}
The restrictive case is $k-i=k$; if $q_2\leq \constA q^{k/(2k-1)}$, then
for $q$ large enough, \eqref{eqn:alphaRootBound} and thus
\eqref{eqn:packCond} holds for all $1\leq k-i\leq k$.

Applying Lemma \ref{lem:MReduct} and noting that $q_2=o(q_1)$, the
left-hand side of the volume condition \eqref{eqn:volCond} is equal
to
\begin{align}
\alpha & \sum_{i=0}^k\sum_{j=0}^i E_i(C)
    \frac{\binom{q_1}{j}\binom{q_2}{i-j}}{t^i}
\left(1+
    \frac{t(i-j)}{\sqrt{2}}\frac{\sqrt{
    \ln(q_2^{1/3}q^{k+1})}}{q_2^{1/3}}(1+o(1))\right).
    \label{eqn:lhsVolCondSumAsymp}
\end{align}
Using the identity $\sum_{j=0}^i \binom{q_1}{j}\binom{q_2}{i-j} =
\binom{q}{i}$, noting that all terms $i<k$ are negligible, bounding
$t(i-j)\leq tk$, and bounding $\ln(q_2^{1/3}q^{k+1})\leq (k+2)\ln
q$, \eqref{eqn:lhsVolCondSumAsymp} is
\begin{equation}
\leq \ \alpha  E_k(C)\frac{\binom{q}{k}}{t^k}
    \left(1+\frac{tk(k+1)\sqrt{k+2}}{\sqrt{2}} \frac{\sqrt{
    \ln q}}{q_2^{1/3}}(1+o(1))\right) = t^{q_2}(1-o(1)).
    \label{eqn:lhsVolCondSumAsymp2}
\end{equation}
By the choice of $\alpha$, \eqref{eqn:volCond} holds for $q$
large enough, so that by Theorem \ref{thm:PaulPackVol} Paul can
win the $(n,q_1,q_2,C)$-game for
\begin{align}
n  \leq &\alpha \left ( 1- q^{-\eta_1}\right) t^{q_1}
 =
 \left \lfloor \frac{t^{q_2+k}}{E_k(C) \binom qk} \left ( 1 -
    \constJ\frac{\sqrt{\ln q}}{q^{1/3}_2} \right) \right \rfloor \left
    (1 - q^{-k-1}\right) t^{q_1},\nonumber
\end{align}
where for $q$ large enough, this last quantity is at least the
right-hand side of \eqref{eqn:PaulOrigWinCond}.

For the pathological variant, let $\constM>t^2k(k+1)\sqrt{k+2}/\sqrt{2}$ and
$\constC>t(t-1)k(k+1)\sqrt{k+2}/\sqrt{2}$ be constants such that
$\constM-t(t-1)k(k+1)\sqrt{k+2}/\sqrt{2}>\constJ$, and
$\constC>\constM-\constJ$.
Let
$$ \alpha'=\left\lceil \constM\frac{t^{q_2+k}}
    {E_k(C)\binom{q}{k}}
    \frac{\sqrt{\ln(q)}}{q_2^{1/3}}\right\rceil.
$$
Following the derivation of \eqref{eqn:lhsPackCondAsymp} and
\eqref{eqn:lhsVolCondSumAsymp2}, for $q$ sufficiently large, the
left-hand side of the volume condition \eqref{eqn:pathVolCond} is
\begin{align}
\geq & \left\lfloor \frac{t^{q_2+k}}{E_k(C)\binom{q}{k}}
    \left(1-\constJ\frac{\sqrt{\ln q}}{q_2^{1/3}}\right)\right\rfloor
     E_k(C)\frac{\binom{q}{k}}{t^k}
    \Big(1-\frac{t(t-1)k(k+1)\sqrt{k+2}}{\sqrt{2}}
    \frac{\sqrt{\ln(q)}}{q_2^{1/3}}\nonumber \\
&    \cdot(1+o(1))\Big)  + \left\lceil \constM\frac{t^{q_2+k}}
    {E_k(C)\binom{q}{k}}
    \frac{\sqrt{\ln(q)}}{q_2^{1/3}}\right\rceil  E_k(C)
    \frac{\binom{q}{k}}{t^k}(1-o(1))
    \geq t^{q_2}, \nonumber
\end{align}
by definition of $\alpha'$, $\constM$, and $\constJ$. Therefore for $q$
sufficiently large, by Theorem \ref{thm:PaulPackVol}, Paul can win
the $(n,q_1,q_2,C)^*$-game for
\begin{eqnarray}
n & \geq &(\alpha+\alpha')(1-q^{-\eta_1})t^{q_1}+t^{q}q^{-\eta_1}
    \nonumber \\
& = & \left( \left\lfloor \frac{t^{q_2+k}}{E_k(C)\binom{q}{k}}
    \left(1-\constJ\frac{\sqrt{\ln q}}{q_2^{1/3}}\right)\right\rfloor
    + \left\lceil \constM\frac{t^{q_2+k}}
    {E_k(C)\binom{q}{k}}
    \frac{\sqrt{\ln(q)}}{q_2^{1/3}}\right\rceil\right)
    (1-q^{-k-1})t^{q_1} \nonumber \\
& & +t^qq^{-k-1}. \nonumber
\end{eqnarray}
By definition of $\constC$, for $q$ sufficiently large, this last
quantity is at most the right-hand side of
\eqref{eqn:PaulWinCond}.
\hfill \qed

{\bfseries{Proof of Theorem \ref{thm:CaroleWinCond}.}}
Fix $\eta_1=\eta_2=\eta=k+1$, and define
$\constN:=t(t-1)k\sqrt{k+2}/\sqrt{2}$. By applying Lemma
\ref{lem:MReduct}, and  noting that
$\ln(q_1^{1/3}q^{\eta_1}),\ln(q_2^{1/3}q^{\eta_2}) \leq(k+2)\ln
q$, the right-hand side of \eqref{eqn:CaroleVolCond} is
\begin{align}
\leq t^q & \left(\frac{E_k(C)}{t^k}  \left(
    \binom{q_1}{k}\left(1 - \constN
    \frac{\sqrt{\ln q}}
    {q_1^{1/3}}(1+o(1))\right) \right. \right. \nonumber \\
& \left.\left.    + \binom{q_2}{k}\left(1 -\constN
    \frac{\sqrt{\ln q}}{q_2^{1/3}}
    (1+o(1))\right)\right)
    \right)^{-1}
+ 2q^{-k-1}t^q. \label{eqn:CaroleAsymp1}
\end{align}
In asymptotically achieving the sphere bound, we cannot have
$\min(q_1,q_2)$ too large, since (for example)
$$
\binom{q_2}{k} = \binom{q}{k}
    \left(\frac{q_2}{q}\right)^k
    \left(1-\frac{q_1}{q_2}\frac{k(k-1)}{2q}(1-o(1))\right).
$$
Assuming for convenience that $q_2=\min{(q_1,q_2)}$, since
$2q^{-k-1}t^q$ is asymptotically negligible, the right-hand
side of \eqref{eqn:CaroleAsymp1} becomes
\begin{align}
t^q & \left(\frac{E_k(C)}{t^k}
    \binom{q}{k}\left(1 -
    \left(\frac{kq_2}{q}+\frac{c_{14}\sqrt{\ln q}}{q_1^{1/3}}
    -\left(\frac{q_2}{q}\right)^k\right)(1+o(1))\right)
    \right)^{-1}. \nonumber
\end{align}
The result for the original variant follows by selecting any
$\constD>k$ (or $\constD>0$ when $k=1$) and $\constE>\constN$
and applying Theorem \ref{thm:CarolePackVol}.
For the pathological variant,  the $i=k$ term dominates the right-hand side of \eqref{eqn:PathCaroleVolCond},
which is asymptotically
\begin{align}
= & t^q(1-q^{-k-1}) \left(\frac{E_k(C)}{t^k} \binom{q}{k}
    \left(1 +\frac{tk\sqrt{k+2}}{\sqrt2}\frac{\sqrt{\ln q}}
    {q^{1/3}}(1+o(1))
    \right) \right)^{-1}. \nonumber
\end{align}
The result follows by selecting any $\constF>tk\sqrt{(k+2)/2}$, noting
that $(1-q^{-k-1})$ is asymptotically negligible, and applying
Theorem \ref{thm:CarolePackVol}. \hfill \qed

\section{Concluding remarks\label{sec:conclusion}}

The first asymptotic term of Theorems
\ref{thm:PaulWinCond}-\ref{thm:CaroleWinCond} is the {\em sphere
bound} for liar games (adaptive codes) over $C$. It arises by
counting the expected number of game response strings for which
$y\in[n]$ survives when Paul's partitions are random, and dividing
into the size of the space $T^q$. Paul's embedding strategy in
Theorem \ref{thm:PaulPackVol} can be viewed as a quasirandom
implementation of this notion.

The most important consequence of Theorem \ref{thm:PaulWinCond},
in the language of coding theory, is the existence of
asymptotically perfect adaptive codes for a wide range of
parameters when the total number of errors (lies) is bounded. The
dominating asymptotic term depends only on the number of lie
strings of maximum length in $C$ and not on their shape.

The generality of the channel led us to make trade-offs for
clarity's sake. For example, the second asymptotic term in
\eqref{eqn:PaulOrigWinCond}-\eqref{eqn:PaulWinCond} could be
reduced to $O(q_2^{-1/3})$ by a more careful embedding of $[n]$ in
$T^{q_1}$, and to $O(q_2^{-1/2})$ by assuming that $C$ is closed
under reordering of lie strings. When $t=2$, the so-called BCH
codes \cite{PHB98} provide a superior bound for $A_t(Q,2R+1)$,
allowing the second batch size to be increased to $q_2=\Theta(q)$
without disturbing the form of the result.  When the suffix
channel $S_C(u')$ of every prefix $u'$ of a length $k$ lie string
$u$ in $C$ is non-degenerate, the original variant bound in
Theorem \ref{thm:CaroleWinCond} improves to $n \geq
t^{q+k}(E_k(C)\binom qk)^{-1}  ( 1 + \mathrm{const.}\cdot
\sqrt{\ln q} ( 1/q_1^{1/3} +1/q_2^{1/3} ))$; this form is superior
when $\min(q_1,q_2)=\omega(\sqrt{\ln q} \, q^{2/3})$. Any channel
such as the binary symmetric, unidirectional, or half-lie channel
that is closed under prefixes has this property, for example. We
are optimistic that Theorem \ref{thm:PaulPackVol} could provide a
basis for understanding the case in which the number of lies grows
to infinity, or for improving bounds on the best known
$k$-error-correcting and radius $k$ covering codes.

\section*{Acknowledgments}

The idea for this paper grew from combining a question posed by
Nathan Linial in 2005: ``does it matter if Carole is using a
$Z$-channel for her lies but Paul doesn't know which one?'',
with our presumption that there should be a liar game
corresponding to adaptive block codes with unidirectional
errors. After completing this work, we were informed of the
extension \cite{ACDV08} of \cite{ACD08}, which is a
specialization of the original liar game case of Theorem
\ref{thm:PaulWinCond}.  We thank Ioana Dumitriu, Joel Spencer,
Alexander Vardy, and Ilya Dumer for several helpful
discussions.  We especially thank the anonymous referee for
carefully reviewing the original manuscript and detecting
several issues.

\end{document}